# A Geometric Interpretation of the Neutrosophic Set, A Generalization of the Intuitionistic Fuzzy Set

Florentin Smarandache, University of New Mexico,
Gallup, NM 87301, USA, E-mail: smarand@unm.edu

*Abstract*:
In this paper we give a geometric interpretation of the Neutrosophic Set using the Neutrosophic Cube. Distinctions between the neutrosophic set and intuitionistic fuzzy set are also presented.

## 1. Introduction:
One first presents the evolution of sets from fuzzy set to neutrosophic set. Then one introduces the neutrosophic components T, I, F which represent the membership, indeterminacy, and non-membership values respectively, where $]^-0, 1^+[$ is the non-standard unit interval, and thus one defines the neutrosophic set.

## 2. Short History:
The *fuzzy set* (FS) was introduced by L. Zadeh in 1965, where each element had a degree of membership.

The *intuitionistic fuzzy set* (IFS) on a universe X was introduced by K. Atanassov in 1983 as a generalization of FS, where besides the degree of membership $\mu_A(x) \in [0,1]$ of each element $x \in X$ to a set A there was considered a degree of non-membership $\nu_A(x) \in [0,1]$, but such that
$$\forall x \in X, \mu_A(x) + \nu_A(x) \leq 1. \tag{2.1}$$
According to Deschrijver & Kerre (2003) the *vague set* defined by Gau and Buehrer (1993) was proven by Bustine & Burillo (1996) to be the same as IFS.

Goguen (1967) defined the *L-fuzzy Set* in X as a mapping $X \to L$ such that $(L^*, \leq_L^*)$ is a complete lattice, where $L^* = \{(x_1, x_2) \in [0,1]^2, x_1 + x_2 \leq 1\}$ and $(x_1, x_2) \leq L^* (y_1, y_2) \Leftrightarrow x_1 \leq y_1$ and $x_2 \geq y_2$. The *interval-valued fuzzy set* (IVFS) apparently first studied by Sambuc (1975), which were called by Deng (1989) *grey sets*, and IFS are specific kinds of L-fuzzy sets.

According to Cornelis et al. (2003), Gehrke et al. (1996) stated that "Many people believe that assigning an exact number to an expert's opinion is too restrictive, and the assignment of an interval of values is more realistic", which is somehow similar with the imprecise probability theory where instead of a crisp probability one has an interval (upper and lower) probabilities as in Walley (1991).

Atanassov (1999) defined the *interval-valued intuitionistic fuzzy set* (IVIFS) on a universe X as an object A such that:
$$A = \{(x, M_A(X), N_A(x)), x \in X\}, \tag{2.2}$$
with $M_A: X \to \text{Int}([0,1])$ and $N_A: X \to \text{Int}([0,1])$ (2.3)
and $\forall x \in X, \sup M_A(x) + \sup N_A(x) \leq 1$. (2.4)

Belnap (1977) defined a four-valued logic, with truth (T), false (F), unknown (U), and contradiction (C). He used a billatice where the four components were inter-related.

In 1995, starting from philosophy (when I fretted to distinguish between *absolute truth* and *relative truth* or between *absolute falsehood* and *relative falsehood* in logics, and respectively between *absolute membership* and *relative membership* or *absolute non-membership* and *relative non-membership* in set theory) I began to use the non-standard analysis. Also, inspired from the sport games (winning, defeating, or tie scores), from votes (pro, contra, null/black votes), from positive/negative/zero numbers, from yes/no/NA, from decision making and control theory (making a decision, not making, or hesitating), from accepted/rejected/pending, etc. and guided by the fact that the law of excluded middle did not work any longer in the modern logics, I combined the non-standard analysis with a tri-component logic/set/probability theory and with philosophy (I was excited by paradoxism in science and arts and letters, as well as by paraconsistency and incompleteness in knowledge). How to deal with all of them at once, is it possible to unity them?

I proposed the term "neutrosophic" because "neutrosophic" etymologically comes from "neutrosophy" [French *neutre* < Latin *neuter*, neutral, and Greek *sophia*, skill/wisdom] which means knowledge of neutral thought, and this third/neutral represents the main distinction between "fuzzy" and "intuitionistic fuzzy" logic/set, i.e. the *included middle* component (Lupasco-Nicolescu's logic in philosophy), i.e. the neutral/indeterminate/unknown part (besides the "truth"/"membership" and "falsehood"/"non-membership" components that both appear in fuzzy logic/set). See the Proceedings of the First International Conference on Neutrosophic Logic, The University of New Mexico, Gallup Campus, 1-3 December 2001, at http://www.gallup.unm.edu/~smarandache/FirstNeutConf.htm.

### 3. Definition of Neutrosophic Set:

Let T, I, F be real standard or non-standard subsets of $]^-0, 1^+[$,
with sup T = t_sup, inf T = t_inf,
sup I = i_sup, inf I = i_inf,
sup F = f_sup, inf F = f_inf,
and n_sup = t_sup+i_sup+f_sup,
n_inf = t_inf+i_inf+f_inf.
T, I, F are called *neutrosophic components*.
Let U be a universe of discourse, and M a set included in U. An element x from U is noted with respect to the set M as x(T, I, F) and belongs to M in the following way:
it is t% true in the set, i% indeterminate (unknown if it is) in the set, and f% false, where t varies in T, i varies in I, f varies in F.

### 4. Neutrosophic Cube as Geometric Interpretation of the Neutrosophic Set:

The most important distinction between IFS and NS is showed in the below **Neutrosophic Cube** A'B'C'D'E'F'G'H' introduced by J. Dezert in 2002.

Because in technical applications only the classical interval $[0,1]$ is used as range for the neutrosophic parameters $t, i, f$, we call the cube $ABCDEDGH$ the **technical neutrosophic cube** and its extension $A'B'C'D'E'D'G'H'$ the **neutrosophic cube** (or **absolute neutrosophic cube**), used in the fields where we need to differentiate between *absolute* and *relative* (as in philosophy) notions.

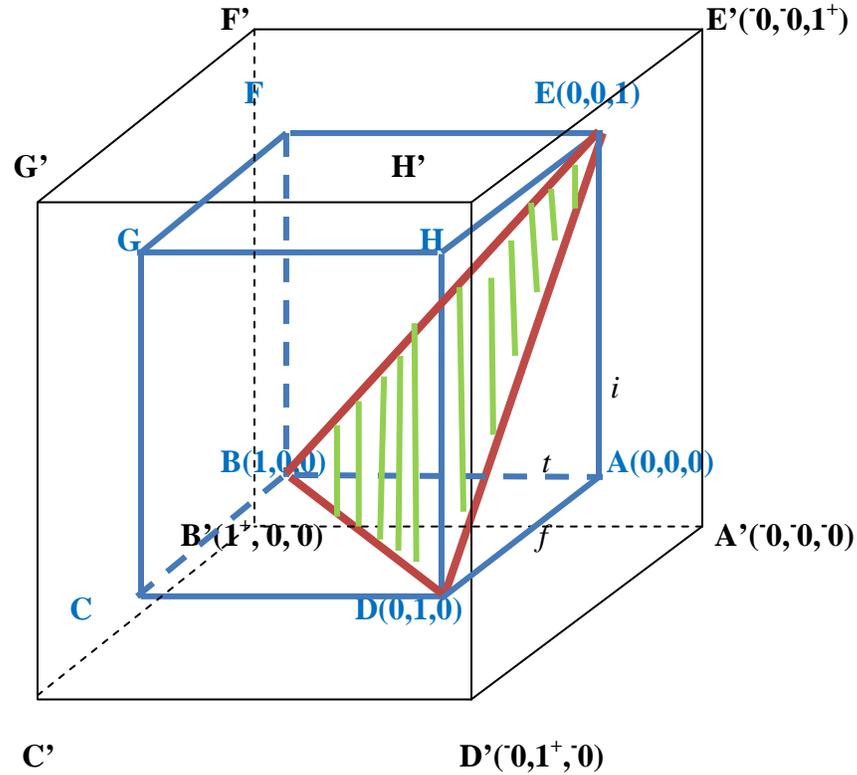

Let's consider a 3D Cartesian system of coordinates, where $t$ is the truth axis with value range in $]^-0,1^+[$, $f$ is the false axis with value range in $]^-0,1^+[$, and similarly $i$ is the indeterminate axis with value range in $]^-0,1^+[$.

We now divide the technical neutrosophic cube $ABCDEDGH$ into three disjoint regions:

1) The equilateral triangle $BDE$, whose sides are equal to $\sqrt{2}$, which represents the geometrical locus of the points whose sum of the coordinates is 1.

If a point $Q$ is situated on the sides of the triangle $BDE$ or inside of it, then $t_Q + i_Q + f_Q = 1$ as in Atanassov-intuitionistic fuzzy set $(A-IFS)$.

2) The pyramid $EABD$ {situated in the right side of the $\Delta EBD$, including its faces $\Delta ABD$ (base), $\Delta EBA$, and $\Delta EDA$ (lateral faces), but excluding its face $\Delta BDE$} is the locus of the points whose sum of coordinates is less than 1.

If $P \in EABD$ then $t_P + i_P + f_P < 1$ as in intuitionistic set (with incomplete information).

3) In the left side of $\Delta BDE$ in the cube there is the solid $EFGCDEBD$ ( excluding $\Delta BDE$ ) which is the locus of points whose sum of their coordinates is greater than 1 as in the paraconsistent set.

If a point $R \in EFGCDEBD$, then $t_R + i_R + f_R > 1$.

It is possible to get the **sum of coordinates strictly less than 1 or strictly greater than 1**. For example:

- We have a source which is capable to find only the degree of membership of an element; but it is unable to find the degree of non-membership;
- Another source which is capable to find only the degree of non-membership of an element;
- Or a source which only computes the indeterminacy.

Thus, when we put the results together of these sources, it is possible that their sum is not 1, but smaller or greater.

Also, in information fusion, when dealing with indeterminate models (i.e. elements of the fusion space which are indeterminate/unknown, such as intersections we don't know if they are empty or not since we don't have enough information, similarly for complements of indeterminate elements, etc.): if we compute the believe in that element (truth), the disbelieve in that element (falsehood), and the indeterminacy part of that element, then the sum of these three components is strictly less than 1 (the difference to 1 is the missing information).

## 5. More distinctions between the Neutrosophic Set and Intuitionistic Fuzzy Set

a) Neutrosophic Set can distinguish between *absolute membership* (i.e. membership in all possible worlds; we have extended Leibniz's absolute truth to absolute membership) and *relative membership* (membership in at least one world but not in all), because NS(absolute membership element)=$1^+$ while NS(relative membership element)=1. This has application in philosophy (see the neutrosophy). That's why the unitary standard interval [0, 1] used in IFS has been extended to the unitary non-standard interval $]^-0, 1^+[$ in NS.
Similar distinctions for *absolute or relative non-membership*, and *absolute or relative indeterminant appurtenance* are allowed in NS.

b) In NS there is no restriction on T, I, F other than they are subsets of $]^-0, 1^+[$,
thus: $^-0 \leq \inf T + \inf I + \inf F \leq \sup T + \sup I + \sup F \leq 3^+$.
The inequalities (2.1) and (2.4) of IFS are relaxed in NS.
This non-restriction allows paraconsistent, dialetheist, and incomplete information to be characterized in NS – as in above Neutrosophic Cube - {i.e. the sum of all three components if they are defined as points, or sum of superior limits of all three components if they are defined as subsets can be >1 (for paraconsistent information coming from different sources), or < 1 for incomplete information}, while that information cannot be described in IFS because in IFS the components T (membership), I (indeterminacy), F (non-membership) are restricted either to t+i+f=1 or to $t^2 + f^2 \leq 1$, if T, I, F are all reduced to the points t, i, f respectively, or to sup T + sup I + sup F = 1 if T, I, F are subsets of [0, 1].
Of course, there are cases when paraconsistent and incomplete informations can be normalized to 1, but this procedure is not always suitable.

c) Relation (2.3) from interval-valued intuitionistic fuzzy set is relaxed in NS, i.e. the intervals

do not necessarily belong to Int[0,1] but to [0,1], even more general to ]⁻0, 1⁺[.

d) In NS the components T, I, F can also be **non-standard subsets** included in the unitary non-standard interval ]⁻0, 1⁺[, not only *standard* subsets included in the unitary standard interval [0, 1] as in IFS.

e) NS, like dialetheism, can describe **paradoxist elements**, NS(paradoxist element) = (1, I, 1), while IFL cannot describe a paradox because the sum of components should be 1 in IFS.

f) The connectors in IFS are defined with respect to T and F, i.e. membership and nonmembership only (hence the Indeterminacy is what's left from 1), while in NS they can be defined with respect to any of them (no restriction).

g) Component "I", indeterminacy, can be split into more subcomponents in order to better catch the vague information we work with, and such, for example, one can get more accurate answers to the *Question-Answering Systems* initiated by Zadeh (2003).

{In Belnap's four-valued logic (1977) indeterminacy is split into Uncertainty (U) and Contradiction (C), but they were interrelated.}

Even more, one can split "I" into Contradiction, Uncertainty, and Unknown, and we get a five-valued logic.

In a general **Refined Neutrosophic Set**, "T" can be split into subcomponents $T_1, T_2, ..., T_m$, and "I" into $I_1, I_2, ..., I_n$, and "F" into $F_1, F_2, ..., F_p$.

h) NS has a better and clear terminology (name) as "neutrosophic" (which means the neutral part: i.e. neither true/membership nor false/nonmembership), while IFS's name "intuitionistic" produces confusion with Intuitionistic Logic, which is something different (see the article by Didier Dubois et al., 2005).

i) The **Neutrosophic Numbers** have been introduced by W.B. Vasantha Kandasamy and F. Smarandache, which are numbers of the form N = a+bI, where a, b are real or complex numbers, while "I" is the indeterminacy part of the neutrosophic number N, such that $I^2 = I$ and $\alpha I + \beta I = (\alpha + \beta)I$.

Of course, indeterminacy "I" is different from the imaginary $i = \sqrt{-1}$.
In general one has $I^n = I$ if n > 0, and is undefined if n ≤ 0.
The algebraic structures using neutrosophic numbers gave birth to the **neutrosophic algebraic structures** [see for example "neutrosophic groups", "neutrosophic rings", "neutrosophic vector space", "neutrosophic matrices, bimatrices, …, n-matrices", etc.], introduced by W.B. Vasantha Kandasamy, F. Smarandache et al.

Example of Neutrosophic Matrix: $\begin{bmatrix} 1 & 2+I & -5 \\ 0 & 1/3 & I \\ -1+4I & 6 & 5I \end{bmatrix}$.

Example of Neutrosophic Ring: ({a+bI, with a, b ∈ R}, +, ·), where of course (a+bI)+(c+dI) = (a+c)+(b+d)I, and (a+bI) · (c+dI) = (ac) + (ad+bc+bd)I.

j) Also, "I" led to the definition of the **neutrosophic graphs** (graphs which have at least either one indeterminate edge or one indeterminate node), and **neutrosophic trees** (trees which have at least either one indeterminate edge or one indeterminate node), which have many applications in social sciences.
As a consequence, the neutrosophic cognitive maps and neutrosophic relational maps are generalizations of fuzzy cognitive maps and respectively fuzzy relational maps (W.B. Vasantha Kandasamy, F. Smarandache et al.).
A Neutrosophic Cognitive Map is a neutrosophic directed graph with concepts like policies, events etc. as nodes and causalities or indeterminates as edges. It represents the causal relationship between concepts.

An edge is said indeterminate if we don't know if it is any relationship between the nodes it connects, or for a directed graph we don't know if it is a directly or inversely proportional relationship.
A node is indeterminate if we don't know what kind of node it is since we have incomplete information.

Example of Neutrosophic Graph (edges $V_1V_3$, $V_1V_5$, $V_2V_3$ are indeterminate and they are drawn as dotted):

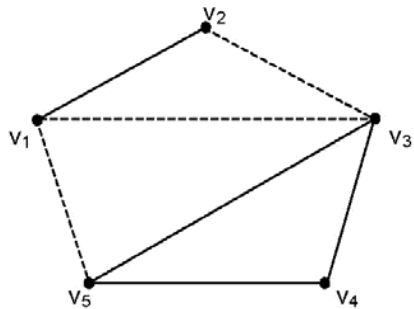

and its neutrosophic adjacency matrix is:

$$\begin{bmatrix} 0 & 1 & I & 0 & I \\ 1 & 0 & I & 0 & 0 \\ I & I & 0 & 1 & 1 \\ 0 & 0 & 1 & 0 & 1 \\ I & 0 & 1 & 1 & 0 \end{bmatrix}$$

The edges mean: 0 = no connection between nodes, 1 = connection between nodes, I = indeterminate connection (not known if it is or if it is not).
Such notions are not used in the fuzzy theory.

Example of Neutrosophic Cognitive Map (NCM), which is a generalization of the Fuzzy Cognitive Maps.
Let's have the following nodes:
$C_1$ - Child Labor
$C_2$ - Political Leaders
$C_3$ - Good Teachers
$C_4$ - Poverty
$C_5$ - Industrialists
$C_6$ - Public practicing/encouraging Child Labor
$C_7$ - Good Non-Governmental Organizations (NGOs)

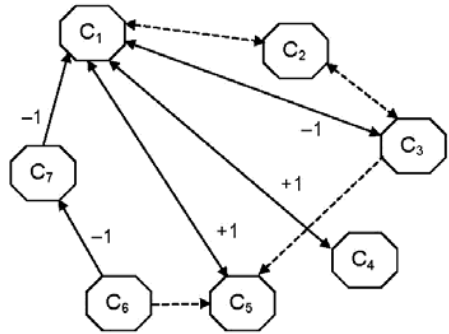

The corresponding neutrosophic adjacency matrix related to this neutrosophic cognitive map is:

$$\begin{bmatrix} 0 & I & -1 & 1 & 1 & 0 & 0 \\ I & 0 & I & 0 & 0 & 0 & 0 \\ -1 & I & 0 & 0 & I & 0 & 0 \\ 1 & 0 & 0 & 0 & 0 & 0 & 0 \\ 1 & 0 & 0 & 0 & 0 & 0 & 0 \\ 0 & 0 & 0 & 0 & I & 0 & -1 \\ -1 & 0 & 0 & 0 & 0 & 0 & 0 \end{bmatrix}$$

The edges mean: 0 = no connection between nodes, 1 = directly proportional connection, -1 = inversely proportionally connection, and I = indeterminate connection (not knowing what kind of relationship is between the nodes the edge connects).

k) The neutrosophics introduced (in 1995) the **Neutrosophic Probability** (NP), which is a generalization of the classical and imprecise probabilities. NP of an event $\mathscr{E}$ is the chance that event $\mathscr{E}$ occurs, the chance that event $\mathscr{E}$ doesn't occur, and the chance of indeterminacy (not knowing if the event $\mathscr{E}$ occurs or not).

In classical probability $n_{sup} \leq 1$, while in neutrosophic probability $n_{sup} \leq 3^+$.

In imprecise probability: the probability of an event is a subset T in [0, 1], not a number p in [0, 1], what's left is supposed to be the opposite, subset F (also from the unit interval [0, 1]); there is no indeterminate subset I in imprecise probability.

And consequently the **Neutrosophic Statistics**, which is the analysis of the neutrosophic events. Neutrosophic statistics deals with neutrosophic numbers, neutrosophic probability distribution, neutrosophic estimation, neutrosophic regression.

The function that models the neutrosophic probability of a random variable x is called *neutrosophic distribution*: $NP(x) = (T(x), I(x), F(x))$, where $T(x)$ represents the probability that value x occurs, $F(x)$ represents the probability that value x does not occur, and $I(x)$ represents the indeterminate / unknown probability of value x.

l) Neutrosophy opened a new field in philosophy.

Neutrosophy is a new branch of philosophy that studies the origin, nature, and scope of neutralities, as well as their interactions with different ideational spectra.

This theory considers every notion or idea <A> together with its opposite or negation <Anti-A> and the spectrum of "neutralities" <Neut-A> (i.e. notions or ideas located between the two extremes, supporting neither <A> nor <Anti-A>). The <Neut-A> and <Anti-A> ideas together are referred to as <Non-A>.

According to this theory every idea <A> tends to be neutralized and balanced by <Anti-A> and <Non-A> ideas - as a state of equilibrium.

In a classical way <A>, <Neut-A>, <Anti-A> are disjoint two by two.

But, since in many cases the borders between notions are vague, imprecise, Sorites, it is possible that <A>, <Neut-A>, <Anti-A> (and <Non-A> of course) have common parts two by two as well.

Neutrosophy is the base of neutrosophic logic, neutrosophic set, neutrosophic probability and statistics used in engineering applications (especially for software and information fusion), medicine, military, cybernetics, physics.

**References on Neutrosophics**